\documentclass{amsart}
\usepackage{amssymb,latexsym} 

\begin{document}
\title{\bf The denominators of harmonic numbers}
\author{Peter Shiu}
\address{353 Fulwood Road, Sheffield, S10 3BQ, United Kingdom}
\email{p.shiu@yahoo.co.uk}
\date{29 July 2024}\begin{abstract} 
The denominators $d_n$ of the harmonic number $1+\frac12+\frac13+\cdots+\frac1n$ do not increase monotonically with~$n$. 
It is conjectured that $d_n=D_n={\rm LCM}(1,2,\ldots,n)$ infinitely often. 
For an odd prime $p$, the set $\{n:pd_n|D_n\}$ has a harmonic density. 
Moreover, for $2<p_1<p_2<\cdots<p_k$, with $\log p_1/\log p_i$ ($1\le i\le k$) being linearly independent, 
there exists $n$ such that $p_1p_2\cdots p_kd_n|D_n$. 
 \end{abstract}
\maketitle

\def\ee#1{\mathcal{E}_#1}
\def\qq#1{\mathcal{Q}_#1}
\def\nd#1#2{\hbox{$#1\!\!\not|#2$}}
\def\nds#1#2{\hbox{$\scriptscriptstyle#1\not\kern1.5pt|#2$}}
\parskip=3pt plus 1pt  minus0.3pt
\section{Introduction}\label{intro}
Although much is known concerning the asymptotic behaviour of the harmonic number
\[H_n=1+\frac12+\frac13+\cdots+\frac1n,\]
there is a dearth of results on the number itself as a fraction.
Let $c_n,d_n$ and~$D_n$ be defined by
\[H_n=\frac{c_n}{d_n},\quad{\rm GCD\,}(c_n,d_n)=1;\qquad D_n={\rm LCM}(1,2,\ldots,n);\]
the following is a small table for their values: 
\begin{center}
\begin{tabular}{|| c | c | c | c | c | c | c | c | c | c | c  ||}\hline
$n$ & $1$ & $2$ & $3$ & $4$ & $5$ & $6$ & $7$ & $8$ & $9$ & $10$ \\ \hline
$c_n$ & $1$ & $3$ & $11$ & $25$ & $137$ & $49$ & $363$ & $761$ & $7129$ & $7381$ \\ \hline
$d_n$ & $1$ & $2$ & $6$ & $12$ & $60$ & $20$ & $140$ & $280$ & $2520$ & $2520$ \\ \hline
$D_n$ & $1$ & $2$ & $6$ & $12$ & $60$ & $60$ & $420$ & $840$ & $2520$ & $2520$ \\ \hline
\end{tabular}
\end{center}
\centerline{\tt Table 1}

Wolstenholme's theorem \cite[Theorem~115]{HandW} states that $p^2|c_{p-1}$ for primes $p>3$. 
More recently, A.~Eswarathasan and E.~Levine~\cite{EandL} conjectured that each odd prime $p$ divides only finitely many~$c_n$, and D.~W.~Boyd~\cite{Boyd} has given computation results and heuristic arguments based on $p$-adic analysis to support the conjecture. 
If the conjecture is true, then there exists $N=N(p)$ such that \nd{p}{c_m} for all $m\ge N$; now, for $n\ge pN$, we may set $n=mp+r$, $0\le r<p$, so that there are precisely $m$ terms in the sum for~$H_n$ in which the denominator is a multiple of~$p$, with their sum being $H_m/p$, and the denominator of the reduced fraction for $H_n-H_m/p$ is free of the prime $p$. It then follows from \nd{p}{c_m} that $p|d_n$ for $n\ge pN$. \par

The sequences $(c_n)$ and $(d_n)$ are not monotonic, and we  prove the following:\par

\goodbreak
{\bf Theorem 1.} \emph{For all $n>1$, we have $c_n\ne c_{n-1}$. Also, each of the following holds for infinitely many~$n$}:\par
(i) $d_n>d_{n-1}$,\  (ii) $d_n=d_{n-1}$,\ (iii) $d_n<d_{n-1}$; \quad\rm (iv) $c_n>c_{n-1}$,\ \rm (v) $c_n<c_{n-1}$.

\smallskip
It is easy to show that the exponent of 2 in the prime factorisation of $d_n$ is the same as that in~$D_n$. Write 
\[D_n=d_nq_n,\]
so that $c_n$ and $q_n$ are odd numbers; and, for odd primes $p$, define the sets  
\[\ee p=\{n: 1<n<p,\,p|c_n\},\qquad\qq p=\{n:p|q_n\}.\] 
Summing the identity 
\begin{equation}\label{symm}\frac1j+\frac1{p-j}=\frac{p}{j(p-j)}, \qquad 1\le j\le\frac{p-1}2,\end{equation}
over $j$, the left-hand side delivers $H_{p-1}$, so that $p|c_{p-1}$, which is enough for our purpose without appealing to Wolstenholme's theorem. 
Anyway,  $p-1\in\ee p$, and we have the following theorems concerning~$\qq p$:

\smallskip
{\bf Theorem 2}. \emph{Let $m\in\ee p$. Then $n\in\qq p$ for each $n$ satisfying
\begin{equation}\label{interval}
mp^a\le n<(m+1)p^a,\qquad a=1,2,\ldots\,.\end{equation}
Conversely, if $n\in\qq p$ then {\rm(\ref{interval})} holds for some $m\in\ee p$ and $a\ge1$.}

\smallskip
{\bf Theorem 3.} \emph{The set $\qq p$ has the harmonic density 
\[\delta(\qq p)=\frac{1}{\log p}\sum_{m\in\ee p}\log\Big(1+\frac1m\Big). \] 
}

\vskip-10pt
{\bf Theorem 4.} \emph{Let $2<p_1<p_2<\cdots<p_k$ be prime numbers. Suppose that 
\[\theta_i=\frac{\log p_1}{\log p_i},\quad i=1,2,\ldots,k, \]
are linearly independent. 
Then there exists $q_n$ such that $p_1p_2\cdots p_k|q_n$.}

The hypothesis that $\theta_i$ being linearly independent is probably unnecessary---it is a consequence of Schanuel's conjecture (see~\cite[pp. 30--31]{Lang}).

\medskip
Since $q_n=1$ means that $n\not\in\qq p$ for all $p\ge3$, we  also define
\[\tilde{\mathcal{Q}}=\text{Comp}\,\bigcup_{p\ge3}\qq p,\qquad
\tilde{Q}(x)=\sum_{n\le x,\,n\in\tilde{\mathcal{Q}}}1,\]
so that $d_n=D_n$ if and only if $n\in\tilde{\mathcal{Q}}$. 
It seems likely that there are infinitely many prime powers $p^a\in \tilde{\mathcal{Q}}$, and we propose the following: 

{\bf Conjecture}. \emph{There are positive constants $K_1,K_2$ such that
\[\frac{K_1x}{\log x}<\tilde{Q}(x)<\frac{K_2x}{\log x},\qquad x>1.\]}

\section{Proof of Theorem 1}\label{pth1} 
(i) The last term in the sum for~$H_p=c_p/d_p$, namely $1/p$, is the only term in which the denominator is a multiple of~$p$.  It follows that $p|d_p$, so that $d_p\ge p$ is unbounded, and hence $d_{n}>d_{n-1}$ for infinitely many~$n$. \par

(ii) We prove that $d_n=d_{n-1}$ when $n=2p>6$. 
Again, $p|d_{n-1}$ because $1/p$ is the only term in the sum for $H_{n-1}$ in which the denominator is a multiple of~$p$. 
Next, the only such terms in the sum for $H_{n}$ are $1/p$ and~$1/2p$, with their sum being $3/p$;  since $p>3$, we deduce that $p|d_n$. Indeed, with $d_{n-1}$ and $d_{n}$ being even, we further deduce that $n=2p$ divides both of them. Finally, it follows from  
\[\frac{c_n}{d_n}=\frac{c_{n-1}+d_{n-1}/n}{d_{n-1}}\qquad\text{and}\qquad
\frac{c_{n-1}}{d_{n-1}}=\frac{c_{n}-d_{n}/n}{d_{n}}\]
that $d_n|d_{n-1}$ and $d_{n-1}|d_n$, so that $d_{n}=d_{n-1}$. \par

(iii) We prove that $d_n<d_{n-1}$ when $n=p(p-1)$. There are precisely $p-2$ terms in the sum for $H_{n-1}$ in which the denominators are multiples of~$p$; in fact
\begin{equation*}
H_{n-1}=\frac{H_{p-2}}p+S_0+S_1+\cdots+S_{p-2},\end{equation*}
where 
\[S_i=\frac1{ip+1}+\frac1{ip+2}+\cdots+\frac1{ip+p-1},\qquad i=0,1,2,\ldots,p-2,\] 
and an obvious generalisation of~(\ref{symm}) shows that
the numerator of the reduced fraction for each $S_i$ is a multiple of~$p$. From  
\[H_{n-1}=\frac{c_{n-1}}{d_{n-1}},\qquad
\frac{H_{p-2}}p=\frac1p\Big(H_{p-1}-\frac1{p-1}\Big)=\frac{c_{p-1}}{pd_{p-1}}-\frac1n,\]
we now find that 
\[\frac{c_n}{d_n}=\frac{c_{n-1}}{d_{n-1}}+\frac1n=\frac{c_{p-1}}{pd_{p-1}}+S_0+S_1+\cdots+S_{p-2}.\]
Since $p|c_{p-1}$, the denominator of the reduced fraction on the right-hand side of the equation is free of the prime~$p$, so that \nd{p}{d_n}, and also $p|d_{n-1}$, because $p|n$. 
Finally, from $d_n|nd_{n-1}$, we now have $p^2d_n|nd_{n-1}$, so that $d_n\le nd_{n-1}/p^2<d_{n-1}$. 

We do not require $p|c_n$ when $p>3$, which follows from Wolstenholme's theorem.  

(iv) If $d_n\ge d_{n-1}$ then $c_n=H_nd_n>H_{n-1}d_{n-1}=c_{n-1}$. 

(v) If $n=p(p-1)$ then $d_n=d_{n-1}/p$. From 
\[\frac{c_{n-1}}{d_{n-1}}=\frac{c_n}{d_n}-\frac1n,\quad\text{that is}\quad 
\frac{c_{n-1}}{c_n}=\frac{d_{n-1}}{d_n}-\frac{d_{n-1}}{nc_n}=p-\frac{d_{n-1}}{nc_n},\]
and
\[\frac{d_{n-1}}{nc_n} =\frac{d_n}{nc_n}{\cdot}\frac{d_{n-1}}{d_n}=\frac{p}{nH_n},\quad\text{so that}\quad
\frac{c_{n-1}}{c_n}=p\Big(1-\frac{1}{nH_n}\Big)>1,\]
we deduce that $c_{n-1}>c_n$ for such $n$. \par

Finally, if $c_n=c_{n-1}$ for any $n>1$, then 
\[d_n=\frac{c_n}{H_n}=\frac{c_{n-1}}{H_n}=\frac{c_{n-1}}{H_{n-1}} {\cdot} \frac{H_{n-1}}{H_n}<d_{n-1} ,\]
which implies that $n$ is composite and that the largest squarefree divisor of $d_{n-1}$ is at least that of~$d_n$. It then follows that $d_{n-1}\ge3d_n$, and hence
\[H_{n-1}+\frac1n=H_n=\frac{c_n}{d_n}=\frac{c_{n-1}}{d_n}\ge\frac{3c_{n-1}}{d_{n-1}}=3H_{n-1},\]
which is impossible. Therefore $c_n\ne c_{n-1}$ for all $n$. \qed

\section{Proof of Theorem 2}\label{pth2} 
Take any $m\in\ee p$, and let $n$ satisfy~(\ref{interval}). 
Then $n>p^a$ so that $p^a|D_n$, and we have
\[H_n=\frac{H_m}{p^a}+\sum_{k\le n\atop\nds{p^a\!}{\,k}}\frac1k
=\frac{c_m/p}{p^{a-1} d_m}+\sum_{k\le n\atop\nds{p^a\!}{\,k}}\frac1k.\]
From $p|c_m$, and hence \nd{p}{d_m}, it follows that \nd{p^a}{d_n}, so that $p|q_n$, that is $n\in\qq p$. \par

Conversely, let $n\in\qq p$, and let $a$ and $m$ be defined by $p^a\le n<p^{a+1}$ and $m=\lfloor n/p^a\rfloor$, so that $1\le m<p$, and $n$ satisfies~(\ref{interval}). 
If $m\not\in\ee p$ then $p^a|d_n$, and since \nd{p^{a+1}}{D_n}, it follows that \nd{p}{q_n}, contradicting $n\in\qq p$. \qed

\goodbreak

\section{Proof of Theorem 3}\label{pth3} 
The length of the interval~(\ref{interval}) is~$p^a$ so that, for each $m\in\ee p$, the number of $n\in\qq p$ with $n\le x=p^b$ is $p+p^2+\cdots+p^{b-1}=(p^b-p)/(p-1)$. 
The intervals corresponding to different $m\in\ee p$ are disjoint, so that 
\begin{equation}\label{asymq}
\qq p(x)=\frac{|\ee p|(p^b-p)}{p-1}\sim\frac{|\ee p|\,x}p\qquad\text{as}\quad x=p^b\to\infty.
\end{equation}
However, since members of $\qq p$ lie in long consecutive runs of integers, the set does not possess an asymptotic density.  
On the other hand its harmonic density exists. As $x\to\infty$, we have
\[\sum_{x\le n<y}\frac1n=\log\frac yx+O\Big(\frac1x\Big),\]
so that
\[\sum_{mp^a\le n<(m+1)p^a}\frac1n=\log\frac{m+1}m+O\Big(\frac1{mp^a}\Big),\]
and hence
\[\sum_{m\in\ee p}\sum_{1\le a<b}\sum_{mp^a\le n<(m+1)p^a}\frac1n=b
\sum_{m\in\ee p}\log\Big(1+\frac1m\Big)
+O(\log p).\]
For $x>1$, we choose $b$ so that $p^b\le x<p^{b+1}$, and 
\[\sum_{x<n<p^{b+1}\atop n\in\qq p}\frac1n
\le \log\frac{p^{b+1}}x+O\Big(\frac1x\Big)\le\log p+O\Big(\frac1x\Big).\]
It then follows that the harmonic density for $\qq p$ is given by
\[\delta(\qq p)=\lim_{x\to\infty}\frac1{\log x}\sum_{n\le x\atop n\in\qq p}\frac1n=
\lim_{b\to\infty}\frac1{b\log p}\sum_{n<p^{b+1}\atop n\in\qq p}\frac1n
=\frac{1}{\log p}\sum_{m\in\ee p}\log\Big(1+\frac1m\Big). \qed\]

\section{Proof of Theorem 4}\label{pth4} 
Take $m_i=p_i-1\in\ee{{p_i}}$, $i=1,2,\ldots,k$ in Theorem~2. 
If $n$ belongs to the intersection of the intervals~(\ref{interval}) corresponding to~$m_i$, then $q_n$ is divisible by each~$p_i$.
The idea then is to find a suitable set of integers $a_i$ so that the intervals can be aligned to deliver a non-empty intersection, or better still, to form a nested sequence of intervals. 
That such an alignment is possible follows from an application of Kronecker's theorem \cite[Theorem~443]{HandW} on simultaneous approximations, assuming that $\theta_i$ are linearly independent. \par

{\bf Lemma 1}. {\it Let $0<\delta<1$, and $2<p_1<p_2<\cdots<p_k$ be the primes in the statement of the theorem. 
Then there are integers $a_1>a_2>\cdots>a_k>1$ such that 
\[(1-\delta)p_1^{a_1}<p_{i+1}^{a_{i+1}}<p_i^{a_i}\le p_1^{a_1},\qquad i=1,2,\ldots,k-1.\]}

\emph{Proof}.
By Kronecker's theorem, corresponding to any $\epsilon>0$, there are positive integers $a_1,a_2,\ldots,a_k$ such that
\[\frac{(i-1)\log(1+\epsilon)}{k\log p_i}\le a_1\theta_i-a_i<\frac{i\log(1+\epsilon)}{k\log p_i},\]
so that
\[\frac{(i-1)\log(1+\epsilon)}{k}\le a_1\log p_1-a_i\log p_i<\frac{i\log(1+\epsilon)}{k},\]
and hence
\[a_1\log p_1-\frac{i\log(1+\epsilon)}{k}<a_i\log p_i\le a_1\log p_1-\frac{(i-1)\log(1+\epsilon)}{k}.\]
It follows that the sequence $(a_i\log p_i)$ is decreasing for $i=1,2,\ldots,k$, and bounded by
$a_1\log p_1-\log(1+\epsilon)$ and $a_1\log p_1$, that is
\[\frac{p_1^{a_1}}{1+\epsilon}< p_{i+1}^{a_{i+1}}<p_i^{a_i}\le p_1^{a_1},\qquad i=1,2,\ldots,k-1.\]
That $a_{i+1}<a_i$ is a consequence of $p_{i+1}>p_i$. 
The required result follows by setting $\epsilon=\delta/(1-\delta)$, so that $(1+\epsilon)(1-\delta)=1$. \qed

\emph{Proof of Theorem 4}.  
Replace $a_i$ by $a_i+1$ in the lemma so that
\[(1-\delta)p_1^{a_1+1}<p_{i+1}^{a_{i+1}+1}<p_i^{a_i+1}\le p_1^{a_1+1},\qquad i=1,2,\ldots,k-1,\]
and hence
\[1-\delta=\frac{(1-\delta)p_1^{a_1+1}}{p_1^{a_1+1}}<\frac{p_{i+1}^{a_{i+1}+1}}{p_i^{a_i+1}}<1.\]
Choosing    $\delta$ to satisfy
\[1-\delta>\min_{1\le i<k}\frac{p_{i+1}(p_i-1)}{p_i(p_{i+1}-1)},\]
we find that
\[\frac{(p_i-1)p_i^{a_i}}{(p_{i+1}-1)p_{i+1}^{a_{i+1}}}
=\frac{p_{i+1}(p_i-1)p_i^{a_i+1}}{p_i(p_{i+1}-1)p_{i+1}^{a_{i+1}+1}}<1,\qquad i=1,2,\ldots,k-1,\]
so that 
\[(p_i-1)p_i^{a_i}<(p_{i+1}-1)p_{i+1}^{a_{i+1}}\le (p_k-1)p_k^{a_k}=(1-1/p_p)p_k^{a_k+1}<p_k^{a_k+1}.\]
With $m_i=p_i-1\in\ee{{p_i}}$, the $k$ intervals 
\[\mathcal{I}(i)=\{n:m_ip_i^{a_i}\le n<(m_i+1)p_i^{a_i}\},\qquad i=1,2,\ldots,k,\]
form a nested sequence, with $\mathcal{I}(i+1)\subset \mathcal{I}(i)$. 
By~(\ref{interval}) in Theorem~2, $p_i|q_n$ for every $n\in\mathcal{I}(i)$, so that $p_1p_2\cdots p_k|q_n$ for every $n\in\mathcal{I}(k)$. The theorem is proved. \qed

\section{Fluctuations in $d_n$}\label{fluc}
The sequence $(D_n)$ is increasing, and $D_{n-1}<D_n$ if and only if $n=p^a$, a prime power---it follows that there are arbitrarily long runs of $D_n$ taking the same value.  
Indeed the inequalities $D_{n-1}<D_n<D_{n+1}$ imply that both $n$ and $n+1$ are prime powers, so that
$\{n,n+1\}=\{2^a,p^b\}$; for example, $n=16$ or~$n=31$.  

The situation for $(d_n)$ is not so simple, and there are interesting problems related to Theorem~1. Thus we may ask if there are arbitrarily long runs of $d_n$ which are strictly increasing,  or  strictly decreasing. 
For example, we find that $d_n<d_{n+1}<d_{n+2}<d_{n+3}$ when $n=1, 2, 6, 22, 70, 820, 856, 1288$. 
On the other hand we have found only one solution to $d_n>d_{n+1}>d_{n+2}$, namely $n=19$; in fact 
$d_{19}=2^43{\cdot}5{\cdot}7{\cdot}11{\cdot}13{\cdot}17{\cdot}19$, and 
$d_{20}=d_{19}/5$, $d_{21}=d_{20}/3$.

A long interval for $n$ in which $D_n$ is stationary may include a long subinterval in which $d_n$ is stationary. For example, $D_n$ takes the same value in $1331\le n<1361$; note that $1331=11^3$, 1361 is a prime, and that $1332=(p-1)p$, with $p=37$. 
The argument in~\S\ref{intro} shows that, for  $n=(p-1)p+r$, $0\le r<p$, the denominator of the reduced fraction for $H_n-H_{p-1}/p$ is free of~$p$, and, since $p|c_{p-1}$, we deduce that \nd{p}{d_n} for $(p-1)p\le n<p^2=1369$; in fact, we find that 
\[D_n=D_{1331}=d_{1331}=37d_n\qquad\text{for}\quad 1332\le n<1361.\]

\section{Listing and counting $n\le x$ with $d_n=D_n$}\label{rrxx}
By the prime number theorem, $D_n=\exp\big((1+o(1))n\big)$ as $n\to\infty$, so that, if $n$ is not small, it is not feasible to compute $H_n$ in order to check whether $d_n=D_n$. 
However, we can apply Theorem~2 to compute~$q_n$ and hence the prime decomposition of~$d_n=D_n/q_n$. 
In particular, we can check whether $n\in\tilde{\mathcal{Q}}$, which amounts to showing that $n$ does not lie in any of the intervals~(\ref{interval}) associated with~$\qq p$. 
If we wish to list $n\in\tilde{\mathcal{Q}}$ in an interval $1\le n\le x$, we can apply a process similar to the sieve of Eratosthenes for the listing of prime numbers.
By this we mean that, for each $m\in\ee p$, each $p\le x/m$ and each $a\le\log x/\log p$, 
we delete, or sift-out, the integers~$n$ in the interval~(\ref{interval}) from the interval \hbox{$1\le n\le x$};  
the un-sifted numbers~$n$ then satisfy \hbox{$d_n=D_n$}.  
Such $n\le 10000$ are given in~\S\ref{results}. \par

It appears that the set $\tilde{\mathcal{Q}}$ can be dealt with using methods applied to the set of primes, and our  Conjecture is based on~(\ref{asymq}).
For example, instead of listing the members of~$\tilde{\mathcal{Q}}$, we may only wish to evaluate, or to estimate, its counting function~$\tilde{Q}(x)$ by applying the inclusion-exclusion principle, in the same way that one does for the  prime counting function~$\pi(x)$. 
The number of $n\le x$ which are divisible by~$a$ has the simple formula~$\lfloor x/a\rfloor$, but it is not so easy to derive a formula for the counting function associated 
with the intersection of various sets~$\qq p$. 
Indeed we have yet to discover why there are arbitrarily large $n$ with $d_n=D_n$; in particular, we have not been able to emulate Euclid's elegant proof that there are infinitely many primes.

\section{Computation results}\label{results} 
There are 2641 numbers $n\le 10000$ such that $d_n=D_n$; they lie in 26 runs of consecutive numbers---in the display below, $a_b$ denotes the interval $a\le n<a+b$.
\begin{align*}
&1_5,\ 9_9,\ 27_6,\ 49_5,\ 88_{12},\ 125_{31},\ 243_{29},\ 289_5,\ 361_2,\ 484_2,\ 841_6,\ 968_{99},\ 1164_2,\\
&1331_1,\ 1369_{89},\ 2401_{99},\ 3125_{254},\ 3488_{172},\ 3721_{272},\ 6889_{87},
\ 7085_{477},\ 7761_{71},\\
&7921_{324},\ 8342_{51},\ 8502_{375},\ 9156_{155}.
\end{align*}
Thus $1_5$ being followed by~$9_9$ means that $1,2,3,4,5,9,10,\ldots,17\in\tilde{\mathcal{Q}}$ and 
\hbox{$6,7,8\not\in\tilde{\mathcal{Q}}$}.  
For $18\le n<27$, there are precisely two terms in the sum for~$H_n$ in which the denominators are divisible by~$3^2$, namely $1/9$ and~$1/18$, with their sum being $1/6$, so that \nd{3^2}{d_n}, whereas~$3^2|D_n$. 
The argument amounts to taking  $p=3$, $m=a=2$ in~(\ref{interval}) to show that such $n\in\qq3$. 
We remark that 15 of the leading terms $a$ for a run are prime powers. 

There is symmetry for~$\ee p$ in that, for $m,m'>0$ with $m+m'+1=p$, we have 
$m\in\ee p$ if and only if $m'\in\ee p$; the proof makes use of the identity~(\ref{symm}).
Since $0\not\in\ee p$ and $p-1\in\ee p$, it follows that $|\ee p|$ is odd, unless $(p-1)/2\in\ee p$. 
A result of Eisenstein~\cite{Eisen} states that
\[-\frac{H_r}2\equiv\frac{2^{p-1}-1}p\pmod p,\qquad r=\frac{p-1}2,\]
so that $(p-1)/2\in\ee p$ if and only if $2^p\equiv2$ (mod~$p^2$).
Such primes $p$ are called Wieferich primes, with only two known ones: $p=1093,3511$, and we find that 
   \[\ee{{1093}}=\{273,546,819,1092\},\qquad \ee{{3511}}=\{877,1755,2633,3510\}.\]
The prime $p<10000$ with the largest  $|\ee p|$ is $p=2113$, and
  \[ \ee{{2113}}=\{44,443,553,748,1384,1559,1669,2068,2112\}.\]
  We end with a table for the number of odd primes $p<10000$ having the same number $|\ee p|$:

\medskip
\begin{center}
\begin{tabular}{|| c | c | c | c | c | c | c | c ||}\hline
\lower6pt\vbox to 18pt{}$|\ee p|$ & $1$ & $3$ & $4$ & $5$ & $7$ & $9$ & Total \\ \hline\hline
\lower6pt\vbox to 18pt{}$|\{3\le p<2039\}|$&192&92& 1&18&4&0&307\\ \hline
\lower6pt\vbox to 18pt{}$|\{2039\le p<4547\}|$&181&97& 1&25&2&1&307\\ \hline
\lower6pt\vbox to 18pt{}$|\{4547\le p<7219\}|$&187&101& 0&17&2&0&307\\ \hline
\lower6pt\vbox to 18pt{}$|\{7219\le p<10000\}|$&192&93& 0&19&3&0&307\\ \hline\hline
\lower6pt\vbox to 18pt{}$|\{3\le p<10000\}|$&752&383& 2&79&11&1&1228\\ \hline
\end{tabular}
\end{center}
\centerline{\tt Table 2}
   
\medskip


\end{document}